\documentclass[11pt]{article}
\usepackage{}
\usepackage{mathrsfs}
\usepackage{amsfonts}

\usepackage{CJK}
\usepackage{amsfonts,amssymb,amsmath,mathrsfs}

\usepackage{color,latexsym,amsfonts}
 \newtheorem{theorem}{Theorem}[section]
 \newtheorem{lemma}[theorem]{Lemma}
 \newtheorem{corollary}[theorem]{Corollary}
 \newtheorem{proposition}[theorem]{Proposition}
 
 \newtheorem{Definition}[theorem]{Definition}
 \newtheorem{remark}[theorem]{Remark}
 \newtheorem{condition}[theorem]{Condition}

 \def\blemma{\begin{lemma}\sl{}\def\elemma{\end{lemma}}}
 \def\btheorem{\begin{theorem}\sl{}\def\etheorem{\end{theorem}}}
 \def\bcorollary{\begin{corollary}\sl{}\def\ecorollary{\end{corollary}}}
 \def\bdefinition{\begin{Definition}\sl{}\def\edefinition{\end{Definition}}}
 \def\bproposition{\begin{proposition}\sl{}\def\eproposition{\end{proposition}}}
 \def\bremark{\begin{remark}\sl{}\def\eremark{\end{remark}}}

 \def\beqlb{\begin{eqnarray}}\def\eeqlb{\end{eqnarray}}
 \def\beqnn{\begin{eqnarray*}}\def\eeqnn{\end{eqnarray*}}

 \def\<{\langle}\def\>{\rangle}

 \def\ar{&\!\!}

 \def\eqref#1{{\rm(\ref{#1})}}

\begin{document}

\

\noindent{}

\bigskip\bigskip

\centerline{\Large\bf A Study of a Class of Stochastic Volterra Equations }

\smallskip

\centerline{\Large\bf Driven by Fractional Brownian Motion}

\smallskip
\

\bigskip\bigskip

\centerline{XiLiang Fan}

\bigskip
\centerline{Department of Mathematics, Anhui Normal University, Wuhu 241003, China}
\smallskip
\centerline{School of Mathematical Sciences, Beijing Normal University, Beijing 100875, China}
\smallskip
\centerline{fanxiliang0515@163.com}

\bigskip\bigskip

{\narrower{\narrower

\noindent{\bf Abstract.} This paper is devoted to study a class of stochastic Volterra equations driven by fractional Brownian motion.
We first prove the Driver type integration by parts formula and the shift Harnack type inequalities.
As a direct application, we provide an alternative method to describe the regularities of the law of the solution.
Secondly, by using the Malliavin calculus, the Bismut type derivative formula is established, which is then
applied to the study of the gradient estimate and the strong Feller property.
Finally, we establish the Talagrand type transportation
cost inequalities for the law of the solution on the path space with respect to both the uniform metric and the $L^2$-metric.
}}


\bigskip



\section{Introduction}

\setcounter{equation}{0}
The Driver integration by parts formula \cite{Driver} and the
Bismut derivative formula \cite{Bismut} are two quite useful tools in various aspects of stochastic analysis.
Let $\nabla$ be the gradient operator and $P_t$ stand for the diffusion semigroup.
The two mentioned formulas allow us to estimate the commutator $\nabla P_t-P_t\nabla$, which plays a key role in the study of flow properties \cite{Fang11}.
On the other hand, \cite{Wang12b} showed that, in general the integration by parts formula is more complicated and harder to obtain than the
derivative formula. Based on martingale method, coupling argument or Malliavin
calculus, the derivative formula has been widely studied and applied
in various fields, such as heat kernel estimates, strong Feller property and functional inequalities,
see \cite{El,Wang10,WangZ,Zhang10a} and references therein. Whereas, in \cite{Wang12b}, based upon
a new coupling argument, the integration by parts formulae are derived and applied to various models
including degenerate diffusion process, delayed SDEs and semi-linear SPDEs.
Afterwards, Zhang (\cite{ZhangSQ12,ZhangSQ13}) studied semi-linear SPDE with delay and stochastic Klein-Gordon type equations; Wang (\cite{Wang13a}) considered SDE with L\'{e}vy noise.

Recently, transportation cost inequality has been widely studied. Let $(E,d)$ be a metric space equipped with $\sigma$-algebra $\mathscr{B}$
such that $d(\cdot,\cdot)$ is $\mathscr{B}\times\mathscr{B}$ measurable. For any $p\geq 1$ and two probability measures $\mu$ and $\nu$ on $(E,\mathscr{B})$,
the $L^p$-Wasserstein distance induced by the metric $d$ between these two probability measures is defined by
\beqnn
W_p^d(\mu,\nu)=\inf\limits_{\pi\in\mathscr{C}(\mu,\nu)}\left(\int_E\int_Ed(x,y)^p\pi({\rm d}x,{\rm d}y)\right)^{1/p},
\eeqnn
where $\mathscr{C}(\mu,\nu)$ denotes the set of all coupling of $\mu$ and $\nu$.
In 1996,  Talagrand \cite{Talagrand} proved the following  transportation cost inequality for the standard Gaussian measure $\mu$ on $\mathbb{R}^d$:
\beqnn
W_2^d(f\mu,\mu)^2\leq2\mu(f\log f), \ f>0, \mu(f)=1,
\eeqnn
where $d(x,y)=|x-y|$.
In general,
we call that the probability measure $\mu$ satisfies the $L^p$-transportation cost inequality on $(E,d)$, if there exists a constant $C(\geq0)$ such that
for any probability measure $\nu$,
\beqlb\label{1.1}
W_p^d(\mu,\nu)\leq\sqrt{2C\mathbb{H}(\nu|\mu)},
\eeqlb
where $\mathbb{H}(\nu|\mu)$, the relative entropy of $\nu$ with respect to $\mu$, is given by
\begin{equation}\label{1.2}\nonumber
\mathbb{H}(\nu|\mu)=
\left\{
\begin{array}{ll}
\int_E\frac{{\rm d}\nu}{{\rm d}\mu}\log\frac{{\rm d}\nu}{{\rm d}\mu}{\rm d}\mu, \ \ {\rm if}\  \nu\ll\mu,\\
+\infty, \ \ \ \ \  \ \ \ \ \ \ \ \ \ {\rm else}.
\end{array} \right.
\end{equation}
For simplicity, we write $\mu\in T_p(C|d)$ for \eqref{1.1}.
In the past decades, the work of Talagrand has been generalized to various  different stochastic processes, see, for instance,
\cite{Otto,Bobkov} for the Hamilton-Jacobi equation, \cite{Feyel02,Feyel04} on abstract Wiener space,
\cite{FangSh} on loop groups, \cite{Djel,WuZh,Wang02,Wang04} for diffusion processes,
 \cite{Wu10} for SDEs of pure jumps, \cite{Ma} for SDEs driven by both Brownian motion and jump process,
\cite{Bao12} for neutral functional SDEs, \cite{Saussereau} for SDEs driven by fractional Brownian motion.

In this article, we are interested in a class of stochastic Volterra equations driven by fractional Brownian motion.
It is well known that the main difficulty raised by the fractional Brownian motion is that it is not Markovian process nor
semimartingale, so the It\^{o} approach to setup a stochastic integral with respect to the fractional Brownian motion is no longer valid.
Now there exist numerous attempts to define a stochastic integral with respect to the fractional Brownian motion and moreover,
many works are devoted to discuss the stochastic differential equations driven by a fractional Brownian motion. We briefly present some results.
Based on a fractional integration by parts formula \cite{Zahle}, Nualart and R\u{a}\c{s}canu \cite{Nualart02} established the existence
and uniqueness result with $H>\frac{1}{2}$. By using the theory of rough path analysis introduced in \cite{Lyons}, Coutin and Qian \cite{Coutin}
proved an existence and uniqueness result with Hurst parameter $H\in(\frac{1}{4},\frac{1}{2})$.
For the regularity results about the law of the solution, the readers may refer to \cite{Hu,Nourdin,Nualart09} and references therein.

The equation associated with a fractional Brownian motion we are to deal with is of Volterra type, which is originally discussed by Coutin and Decreusefond \cite{COutin&Decreu}. In \cite{COutin&Decreu}, the authors studied existence, uniqueness and regularity of solution.
In this paper, by using coupling argument and the Girsanov transform for fractional Brownian motion, we first prove the Driver type integration by parts formula and then derive shift Harnack type inequalities.
As an important application, we give an alternative proof of \cite[Corollary 4.1]{COutin&Decreu}.
Secondly, based on Malliavin calculus the Bismut type derivative formula is established, which is then applied to study the
gradient estimate, the Harnack type inequalities and the strong Feller property.
Finally, we obtain the Talagrand type transportation cost inequalities for the law of the solution
on the path space with respect to both the uniform metric and $L^2$-metric.

The paper is organized as follows. In section 2, we give some preliminaries on fractional Brownian motion.
In section 3, we investigate the Driver type integration by parts formula, while in section 4,
the Bismut type derivative formula is discussed.
Finally, section 5 is devoted to derivation of the transportation cost inequalities.

\section{Preliminaries}

\setcounter{equation}{0}

Let $B^H=\{B_t^H, t\in[0,T]\}$ be a fractional Brownian motion with
Hurst parameter $H\in(0,1)$ defined on a complete probability space
$(\Omega,\mathscr{F},\mathbb{P})$. Namely,  $B^H$ is a centered Gauss
process with the covariance function \beqnn
\mathbb{E}(B_t^{H}B_s^{H})=R_H(t,s):=\frac{1}{2}\left(t^{2H}+s^{2H}-|t-s|^{2H}\right).
\eeqnn
When $H=\frac{1}{2}$, the process $B^\frac{1}{2}$ is the usual Brownian motion.
By the above covariance function and the Kolmogorov criterion,  we know that
$B^{H}$ have $(H-\epsilon)$-order H\"{o}lder continuous paths for all $\epsilon>0$.
Furthermore, $B^{H}$ has stationary increments and is self-similar with respect to Hurst index $H$.

From \cite{Decr}, it is known that the covariance kernel $R_H(t,s)$ admits the following representation:
\beqnn
R_H(t,s)=\int_0^{t\wedge s}K_H(t,r)K_H(s,r){\rm d}r,
\eeqnn
where $K_H(\cdot,\cdot)$ is a square integrable kernel given by
\beqnn
K_H(t,s)=\Gamma\left(H+\frac{1}{2}\right)^{-1}(t-s)^{H-\frac{1}{2}}F\left(H-\frac{1}{2},\frac{1}{2}-H,H+\frac{1}{2},1-\frac{t}{s}\right),\ t>s>0,
\eeqnn
in which $F(\cdot,\cdot,\cdot,\cdot)$ is the Gauss hypergeometric function (for details, see \cite{Niki}).\\
Again by \cite{Decr}, the operator $K_H: L^2([0,T];\mathbb{R})\rightarrow I_{0+}^{H+1/2}(L^2([0,T];\mathbb{R}))$
associated with the kernel $K_H(\cdot,\cdot)$ is defined as follows
\beqnn
(K_Hf)(t):=\int_0^tK_H(t,s)f(s){\rm d}s,
\eeqnn
where $ I_{0+}^{H+1/2}$  is the $(H+1/2)$-order left fractional Riemann-Liouville integral operator on $[0,T]$.
It is an isomorphism and for each $f\in L^2([0,T];\mathbb{R})$,
\beqlb\label{2.1,2.2}
(K_Hf)(s)=I_{0+}^{2H}s^{1/2-H}I_{0+}^{1/2-H}s^{H-\frac{1}{2}}f,\ H\leq1/2,\\
(K_Hf)(s)=I_{0+}^{1}s^{H-1/2}I_{0+}^{H-1/2}s^{1/2-H}f,\ H\geq1/2.
\eeqlb
Hence, for any $h\in I_{0+}^{H+1/2}(L^2([0,T];\mathbb{R}))$, the inverse operator $K_H^{-1}$ can be written as
\beqnn
 (K_H^{-1}h)(s)=s^{H-1/2}D_{0+}^{H-1/2}s^{1/2-H}h',\ H>1/2,\\
 (K_H^{-1}h)(s)=s^{1/2-H}D_{0+}^{1/2-H}s^{H-1/2}D_{0+}^{2H}h,\ H<1/2,
\eeqnn
where $D_{0+}^{\alpha}$ is the $\alpha$-order left-sided Riemann-Liouville derivative operator, $\alpha\in(0,1)$.\\
In particular, when $h$ is absolutely continuous, it holds
\beqnn
 (K_H^{-1}h)(s)=s^{H-1/2}I_{0+}^{1/2-H}s^{1/2-H}h',\ H<1/2.
\eeqnn
For more details about the deterministic fractional calculus, one can refer to \cite{Samko}.

We assume that $\Omega$ is the canonical probability space $C_0([0,T];\mathbb{R})$, the set of continuous functions, null at time 0,
equipped with the Borel $\sigma$-algebra associated with the supremum norm and $\mathbb{P}$ is the law of the fractional Brownian motion. The canonical filtration is
$\mathscr{F}_t=\sigma\{B^H_s:0\leq s\leq t\}\vee\mathcal {N}$, where $\mathcal {N}$ is the set of the $\mathbb{P}$-null sets. According to \cite[Theorem 3.3]{Decr}, the Cameron-Martin
space of the fractional Brownian motion, denoted by $\mathcal {H}$, is equal to $I_{0+}^{H+1/2}(L^2([0,T];\mathbb{R}))$, i.e., for any
$h\in\mathcal {H}$, it can be represented as $h(t)=K_H\dot{h}(t)$, where the function $\dot{h}$ belongs to $L^2([0,T];\mathbb{R})$.
The scalar product on $\mathcal {H}$ is defined by
\beqnn
 (h,g)_\mathcal {H}:=(K_H^{-1}h,K_H^{-1}g)_{L^2([0,T];\mathbb{R})},\  \forall h, g\in\mathcal {H}.
\eeqnn
As a consequence, $(\Omega,\mathcal {H},\mathbb{P})$ is an abstract Wiener space in the sense of Gross.
Furthermore, let $\Omega^\ast$ denote the strong topological dual of $\Omega$, then there hold
\beqnn
\Omega^\ast\xrightarrow{K_H^\ast}L^2([0,T];\mathbb{R})\xrightarrow{K_H}\mathcal {H}\xrightarrow{i_H}\Omega
\eeqnn
and
\beqnn
R_H=K_H\circ K_H^\ast,
\eeqnn
where we identify the operator $R_H$ and its kernel.

Next we summarize some basic results of Malliavin calculus associated with the fractional Brownian motion, and we refer to
\cite{Decr}, \cite{Nualart06} and \cite{Ust} for a comprehensive account.

Let $\mathcal {S}$ denote the set of smooth and cylindrical random variables of the form:
\beqnn
F=f(\langle l_1,\omega\rangle,\cdot\cdot\cdot,\langle l_n,\omega\rangle)
\eeqnn
where $n\geq 1, f\in C_b^\infty(\mathbb{R}^n)$, the set of $f$ and all its partial derivatives are bounded, $l_i\in\Omega^\ast, 1\leq i\leq n$.
The Malliavin derivative of $F$, denoted by $D_HF$, is defined as the $\mathcal {H}$-valued random variable
\beqnn
D_HF(\omega)=\sum_{i=1}^n\frac{\partial f}{\partial x_i}(\langle l_1,\omega\rangle,\cdot\cdot\cdot,\langle l_n,\omega\rangle)R_H(l_i).
\eeqnn
For any $k\in\mathbb{N}$, denote by $D_H^k$ the iteration of $D_H$. For any $p\geq 1$ and $k\in\mathbb{N}$, we define the Sobolev space
$\mathbb{D}_H^{k,p}$ as the completion of $\mathcal {S}$ with respect to the norm:
\beqnn
\|F\|_{k,p}^p:=\mathbb{E}|F|^p+\mathbb{E}\sum_{i=1}^k\|D_H^iF\|^p_{\mathcal {H}^{\otimes i}}.
\eeqnn
The divergence operator $\delta_H$, also called the Skorohod integral, is defined by using the duality relationship. More precisely,
the domain Dom$_p\delta_H$ is the set of process $u$ such that
\beqnn
|\mathbb{E}\langle D_HF, u\rangle_\mathcal{H}|\leq C(\mathbb{E}|F|^q)^\frac{1}{q}
\eeqnn
for all $F\in\mathbb{D}_H^{1,q}$, where $q$ satisfies $1/p+1/q=1$ and $C$ is some constant depending on $u$.\\
If $u\in$ Dom$_p\delta_H$, then $\delta_Hu$ is defined by
\beqnn
\mathbb{E}(F\delta_Hu)=\mathbb{E}\langle D_HF, u\rangle_\mathcal{H}.
\eeqnn

It is well known that, in the case of the Brownian motion $(H=1/2)$,  the Skorohod integral is an extension of the It\^{o} integral.
So, this motivates us to use the divergence operator to define a stochastic integral with respect to the fractional Brownian motion. That is,
\beqnn
\int_0^Tu_t\delta_HB^H_t:=\delta_H(K_Hu),
\eeqnn
where the process $K_Hu\in{\rm Dom}\delta_H:=\cup_{p\geq1}{\rm Dom}_p\delta_H$ (see e.g. \cite{Decr} and \cite{COutin&Decreu}).
According to \cite[Theorem 4.8]{Decr}, we have the following L\'{e}vy-Hida representation:
\beqnn
W:=\left(\int_0^t{\rm I}_{[0,t]}(s)\delta_HB^H_s\right)_{0\leq t\leq T}=(\delta_H(K_H{\rm I}_{[0,t]}))_{0\leq t\leq T}
\eeqnn
is a standard Brownian motion whose filtration is equal to $\{\mathscr{F}_t, 0\leq t\leq T\}$ and moreover,
for any square integrable and adapted process $u$, that is $u\in L_a^2([0,T]\times\Omega)$,
it holds
\beqnn
\int_0^tu_s\delta_HB^H_s=\int_0^tu_s{\rm d}W_s, \forall t\in[0,T].
\eeqnn
In particular, for each $t\in[0,T]$, taking $u=K_H(t,\cdot)$,  then we get $B_t^H=\int_0^tK_H(t,s){\rm d}W_s$.

In present paper, we are concerned with a $\mathbb{R}$-valued equation driven by a fractional Brownian motion of the form:
\beqlb\label{2.3}
 X_t=x+\int_0^tK_H(t,s)b(s,X_s){\rm d}s+\int_0^tK_H(t,s)\sigma(s,X_s){\rm d}W_s.
\eeqlb
Note that, when $\sigma\equiv C$, then the third term of the right-hand side of \eqref{2.3} is equal to $CB^H_t$. Hence, the factor $K_H(t,s)$ in the noise
term is necessary to make the equation well-defined.
While $K_H(t,s)$  in the drift term is only to symmetrize $b$ and $\sigma$.
Set $H_0=|H-1/2|, A_H=\{p\geq1: pH_0<1\}$  and for every  $p\in A_H$, put $\kappa_p=(1-pH_0)^{-1}, L_+^{\kappa_p}=\cup_{q>\kappa_p}L^q([0,T];\mathbb{R})$.
\bdefinition\label{d2.1}
A $\mathbb{R}$-valude process $(X_t)_{t\in[0,T]}$ is called a solution of \eqref{2.3}, if it is adapted such that
$\mathbb{E}|X_t|^2\in L_+^{\kappa_2}$, and \eqref{2.3} is satisfied ${\rm d}\mathbb{P}\times {\rm d}t$ a.s.
\edefinition

\bremark\label{r2.1}
From \cite[Theorem 3.1 and Theorem 3.2]{COutin&Decreu}, we know that, if $b$ and $\sigma$ are Lipschitz continuous for the second variable uniformly with respect to
their first variable, and there exist $x_0, y_0\in\mathbb{R}$ such that
$b(\cdot,x_0)\in L^{\kappa_1}_+([0,T];\mathbb{R}), \sigma(\cdot,y_0)\in L^{2\kappa_2}_+([0,T];\mathbb{R})$, then \eqref{2.3} has a unique solution.
Furthermore, for all $p\in A_H, \sup_{t\in[0,T]}\mathbb{E}|X_t|^p<\infty$. If $\sigma$ is bounded, then $X$ has almost surely continuous trajectories.
\eremark

Define $P_tf(x):=\mathbb{E}f(X_t^x), \ t\in[0,T], \ f\in \mathscr{B}_b(\mathbb{R})$, where $X_t^x$ is the solution to \eqref{2.3}
with $X_0=x$ and $\mathscr{B}_b(\mathbb{R})$ denotes the set of all bounded measurable functions on $\mathbb{R}$.
Besides, we denote by $C_b^1(\mathbb{R})$ the set of all bounded continuous differentiable functions.
In the remainder of the paper, we will establish the Driver type integration by parts formula and the Bismut type derivative formula for $P_T$, and moreover obtain the Talagrand type transportation cost inequalities for the law of the solution of \eqref{2.3} on the path space.

\section{Driver type integration by parts formula}

\setcounter{equation}{0}
This section is devoted to the equation \eqref{2.3} with additive noise, i.e.
\beqlb\label{3.1}
X_t=x+\int_0^tK_H(t,s)b(s,X_s){\rm d}s+\int_0^tK_H(t,s)\sigma(s){\rm d}W_s.
\eeqlb
We aim to establish the Driver type integration by parts formula
and the shift Harnack inequalities by the method of coupling and Girsanov transformation.
As an application, we give an alternative proof of \cite[Corollary 4.1]{COutin&Decreu}, in which the
absolute continuity of the law of the solution is discussed.

To start with, let us introduce the notation
$$\nabla_yf(x):=\lim\limits_{\epsilon\downarrow0}\frac{f(x+\epsilon y)-f(x)}{\epsilon},$$
and give some conditions of the coefficients $b$ and $\sigma$: (H1)
\begin{itemize}
\item[(i)] $b$ is continuously differentiable w.r.t. the second variable and
there exist positive constants $K_1$ and $K_2$ such that
$$|\partial b(t,\cdot)(x)|\leq K_1, |\sigma(t)^{-1}|\leq K_2,\ \forall t\in[0,T], x\in\mathbb{R};$$
  \par
 \item[(ii)] there exist $x_0\in\mathbb{R}$ such that
$b(\cdot,x_0)\in L^{\kappa_1}_+([0,T];\mathbb{R}), \sigma(\cdot)\in L^{2\kappa_2}_+([0,T];\mathbb{R})$.
\end{itemize}

\btheorem\label{t3.1}
Let $T>0$ and $y\in\mathbb{R}$ be fixed. Assume that (H1) holds.
\begin{itemize}
 \item[{\rm(1)}]
For each $f\in C_b^1(\mathbb{R})$, there holds the integration by parts formula
$$P_T(\nabla_yf)(x)=\mathbb{E}\left[f(X_T^x)\int_0^T\sigma(s)^{-1}\left(C_Hs^{\frac{1}{2}-H}-s\partial b(s,\cdot)(X_s)\right)\frac{y}{T}{\rm d}W_s\right],$$
where $C_H$ is a positive constant given in the proof below. \\
 As a consequence, for each $\alpha>0$ and positive $f\in C_b^1(\mathbb{R})$,
 \beqnn
 |P_T(\nabla_yf)|\ar\leq\ar\alpha\left[P_T(f\log f)-(P_Tf)(\log P_Tf)\right]\cr
 \ar+\ar\frac{K_2^2y^2}{\alpha}\left(\frac{C_H^2}{(2-2H)T^{2H}}+\frac{4C_HK_1}{5-2H}T^{\frac{1}{2}-H}+\frac{K_1^2T}{3}\right)P_Tf.
 \eeqnn
  \par
 \item[{\rm(2)}] For each non-negative $f\in\mathscr{B}_b(\mathbb{R})$, there holds the shift Harnack inequality
 $$(P_Tf)^p\leq\left(P_T{\{f(y+\cdot)\}}^p\right)\exp\left[\frac{pK_2^2}{p-1}\left(\frac{C_H^2}{(2-2H)T^{2H}}+\frac{4C_HK_1}{5-2H}T^{\frac{1}{2}-H}+\frac{K_1^2T}{3}\right)y^2\right].$$
 \par
 \item[{\rm(3)}] For each positive $f\in\mathscr{B}_b(\mathbb{R})$, there holds the shift log-Harnack inequality
$$P_T\log f\leq\log P_T\{f(y+\cdot)\}+K_2^2\left(\frac{C_H^2}{(2-2H)T^{2H}}+\frac{4C_HK_1}{5-2H}T^{\frac{1}{2}-H}+\frac{K_1^2T}{3}\right)y^2.$$
\end{itemize}
\etheorem

\emph{Proof.}
Obviously, by (H1), it follows from Remark \ref{r2.1} that \eqref{3.1} has a unique solution.
On the other hand, for any $\epsilon\in[0,1]$, let $X_t^\epsilon$ solve the equation
\beqlb\label{3.2}
X_t^\epsilon=x+\int_0^tK_H(t,s)b(s,X_s){\rm d}s+\int_0^tK_H(t,s)\sigma(s){\rm d}W_s+\frac{t\epsilon}{T}y, \ t\in[0,T].
\eeqlb
It is easy to see that $X_t^\epsilon=X_t+\frac{t\epsilon}{T}y,\ t\in[0,T]$. In particular, $X_T^\epsilon=X_T+\epsilon y$.\\
Next, let
\begin{equation}\nonumber
C_H=
\left\{
\begin{array}{ll}
\frac{\Gamma(2H)\Gamma(1/2-H)}{B(2-2H,2H)}\left(\frac{1}{2}-H\right),\ H<\frac{1}{2};\\
1,\ \ \ \ \ \ \ \ \ \ \ \ \ \ \ \ \ \ \ \ \ \ \ \ \ \ \ \ \ H=\frac{1}{2};\\
\frac{\Gamma(H-1/2)}{B(2-2H,H-1/2)},\ \ \ \ \ \ \ \ \  \ \ \ H>\frac{1}{2}.
\end{array} \right.
\end{equation}
It follows from \eqref{2.1,2.2} and (2.2) that, for each $H\in(0,1), K_H(C_Hx^{1/2-H})(t)=t.$ Therefore, we can reformulate \eqref{3.2} as
\beqnn
X_t^\epsilon=x+\int_0^tK_H(t,s)b(s,X_s^\epsilon){\rm d}s+\int_0^tK_H(t,s)\sigma(s){\rm d}W_s^\epsilon,\ t\in[0,T],
\eeqnn
where
$W_s^\epsilon=W_s+\int_o^s\sigma(r)^{-1}\left(b(r,X_r)-b(r,X_r^\epsilon)+\frac{C_H\epsilon y}{T}r^{\frac{1}{2}-H}\right){\rm d}r,\ s\in[0,T].$\\
Let
\beqnn
\xi_\epsilon(r)=b(r,X_r)-b(r,X_r^\epsilon)+\frac{C_H\epsilon y}{T}r^{\frac{1}{2}-H}
\eeqnn
and
\beqnn
R_\epsilon=\exp\left[-\int_0^T\sigma(r)^{-1}\xi_\epsilon(r){\rm d}W_r-\frac{1}{2}\int_0^T|\sigma(r)^{-1}\xi_\epsilon(r)|^2{\rm d}r\right].
\eeqnn
According to (H1), we easily get $\mathbb{E}\exp\left[\frac{1}{2}\int_0^T|\sigma(r)^{-1}\xi_\epsilon(r)|^2{\rm d}r\right]<\infty.$ By the Novikov condition and the Girsanov theorem, $(W^\epsilon_t)_{0\leq t\leq T}$ is a Brownian motion under the probability measure $\mathbb{Q}_\epsilon:=R_\epsilon\mathbb{P}$. Then $(X,X^\epsilon)$ is a coupling by change of measure with changed probability $\mathbb{Q}_\epsilon$.
Since $R_0=1$, by \cite[Theorem 2.1]{Wang12b}, to obtain the desired integration by parts formula, it remains to confirm the following equality:
in the sense of $L^1(\mathbb{P})$,
\beqnn
\frac{{\rm d}}{{\rm d}\epsilon}R_\epsilon|_{\epsilon=0}=-\int_0^T\sigma(r)^{-1}\left[C_Hr^{\frac{1}{2}-H}-r\partial b(r,\cdot)(X_r)\right]\frac{y}{T}{\rm d}W_r.
\eeqnn
Actually, noting that
\beqnn
\lim\limits_{\epsilon\rightarrow 0}\mathbb{E}\frac{R_\epsilon-1}{\epsilon}
\ar=\ar\lim\limits_{\epsilon\rightarrow 0}\mathbb{E}\frac{-\int_0^T\sigma(r)^{-1}\xi_\epsilon(r){\rm d}W_r-\frac{1}{2}\int_0^T|\sigma(r)^{-1}\xi_\epsilon(r)|^2{\rm d}r}{\epsilon}\cr
\ar=\ar\lim\limits_{\epsilon\rightarrow 0}\mathbb{E}\frac{-\int_0^T\sigma(r)^{-1}\xi_\epsilon(r){\rm d}W_r}{\epsilon},
\eeqnn
and, moreover
\beqnn
\ar\mathbb{E}\ar\left|\frac{-\int_0^T\sigma(r)^{-1}\xi_\epsilon(r){\rm d}W_r}{\epsilon}+\int_0^T\sigma(r)^{-1}\left[C_Hr^{\frac{1}{2}-H}-r\partial b(r,\cdot)(X_r)\right]\frac{y}{T}{\rm d}W_r\right|\cr
\ar\leq\ar K_2\left[\mathbb{E}\int_0^T\left|\frac{b(r,X_r^\epsilon)-b(r,X_r)-\partial b(r,\cdot)(X_r)\frac{ry}{T}\epsilon}{\epsilon}\right|^2{\rm d}r\right]^\frac{1}{2},
\eeqnn
so the dominated convergence theorem implies the assertion.\\
The second result in (1) follows by the given upper bounds on $|\sigma(t)^{-1}|$ and $|\partial b(t,\cdot)|$ and using
the above integration by parts formula and the Young inequality (see, for instance, \cite[Lemma 2.4]{Arnaudon})
\beqnn
\ar\ar|P_T(\nabla_yf)|-\alpha\left[P_T(f\log f)-(P_Tf)(\log P_Tf)\right]\cr
\ar\leq\ar\alpha\log\mathbb{E}\exp\left[\frac{1}{\alpha}\int_0^T\sigma(s)^{-1}\left(C_Hs^{\frac{1}{2}-H}-s\partial b(s,\cdot)(X_s)\right)\frac{y}{T}{\rm d}W_s\right]\cdot P_Tf\cr
\ar\leq\ar\frac{\alpha}{2}\log\mathbb{E}\exp\left[\frac{2}{\alpha^2}\int_0^T\left|\sigma(s)^{-1}\left(C_Hs^{\frac{1}{2}-H}-s\partial b(s,\cdot)(X_s)\right)\frac{y}{T}\right|^2{\rm d}s\right]\cdot P_Tf.
\eeqnn
Finally, (2) and (3) can be easily derived by applying \cite[Propositon 2.3]{Wang12b} and the second inequality in (1). The proof is complete.

The shift Harnack type inequalities allow us to deduce the regularity for the law of the solution of \eqref{3.1}. That is, we have the following result.
\bcorollary\label{c3.1}
Suppose that the assumption (H1) holds. Then, for any $t>0$, the law of $X_t$ is  absolutely continuous with respect to the Lebesgue measure.
\ecorollary

\emph{Proof.}
Without lost of generality, we only consider the case $t=T$. Let
$$C_1=\frac{pK_2^2}{p-1}\left(\frac{C_H^2}{(2-2H)T^{2H}}+\frac{4C_HK_1}{5-2H}T^{\frac{1}{2}-H}+\frac{K_1^2T}{3}\right).$$
The shift Harnack inequality stated in Theorem \ref{t3.1} implies,  for any non-negative $f\in\mathscr{B}_b(\mathbb{R})$,
\beqnn
(P_Tf(x))^p{\rm e}^{-C_1|y|^2}\leq\left(P_T{\{f(y+\cdot)\}}^p\right)(x).
\eeqnn
Let $A$ be a Lebesgue-null set, by applying the above inequality to $f={\rm I}_A$ and noting the invariance property under shift for the Lebesgue measure,
we have
$$(P_T{\rm I}_A(x))^p\int_{\mathbb{R}}{\rm e}^{-C_1|y|^2}{\rm d}y\leq 0,$$
which implies the desired result.

\section{Bismut type derivative formula}

\setcounter{equation}{0}
In this section, we shall adopt the techniques of the Malliavin calculus to investigate the Bismut type derivative formula and the Harnack type inequalities for
$P_T$ associated with \eqref{3.1}.
To this end, we make the following assumption: (H2)
\begin{itemize}
 \item[(i)] there exist $x_0\in\mathbb{R}$ and $p\geq2$ such that $b(\cdot,x_0)\in L^p([0,T];\mathbb{R})$;
\par
 \item[(ii)] $b$ is differentiable w.r.t.
 the space variable such that $\partial b(t,\cdot)$ is uniformly continuous uniformly w.r.t. the time variable $t$
 and moreover,
\beqnn
|\partial b(t,\cdot)(x)|\leq K_3, K_5\leq|\sigma(t)^{-1}|\leq K_4, \ \forall t\in[0,T], x\in\mathbb{R},
\eeqnn
where $K_3, K_4$ and $K_5$ are positive constants.
\end{itemize}

Main result reads as follows.
\btheorem\label{t4.1}
Assume that (H2) holds. Then, for all $x, y\in \mathbb{R}$ and $f\in C_b^1(\mathbb{R})$,
\beqnn
\nabla_yP_Tf(x)=\mathbb{E}\left[f(X_T^x)\int_0^T\sigma(s)^{-1}\left(\left(1+\int_0^sK_H(s,r)u'(r){\rm d}r\right)\partial b(s,\cdot)(X_s^x)-u'(s)\right)y{\rm d}W_s\right],
\eeqnn
where $u\in C^1([0,T];\mathbb{R})$ such that $1+\int_0^TK_H(T,r)u'(r){\rm d}r=0, X_\cdot^x$ is the solution of \eqref{3.1}.
\etheorem

The proof of this theorem is based on the following lemmas and proposition.

We first recall a result from  \cite[Theorem 4.1 and Theorem 4.2]{COutin&Decreu}, in which the existence of Malliavin directional derivative is discussed.
\blemma\label{l4.1}
Let $b$ and $\sigma$ be continuously differentiable w.r.t. their space variable, with bounded derivative;
assume further that, there exist $x_0\in\mathbb{R}$ and $p\geq2$ such that $b(\cdot,x_0)\in L^p([0,T];\mathbb{R})$ and $\sigma$ is bounded. Then, for any
$\xi\in\mathcal {H}, (\langle D_HX_t^x, \xi\rangle_\mathcal {H})_{t\in[0,T]}$ exists and
is the unique solution to the equation
\beqnn
Y_t=\langle K_H(K_H(t,\cdot)\sigma(\cdot,X_\cdot^x)), \xi\rangle_\mathcal {H}
+\int_0^tK_H(t,s)\partial b(s,\cdot)(X_s^x)Y_s{\rm d}s+\int_0^tK_H(t,s)\partial\sigma(s,\cdot)(X_s^x)Y_s{\rm d}W_s,
\eeqnn
where $X_\cdot^x$ is the solution of \eqref{2.3}.
\elemma

Following the same method presented in \cite[Theorem 3.3]{COutin&Decreu}, we can show that the solution of \eqref{2.3} depends continuously on
the initial condition in the sense specified below.
\blemma\label{l4.2}
Assume $b$ and $\sigma$ are Lipschitz continuous for the second variable uniformly w.r.t.
their first variable, and there exist $x_0, y_0\in\mathbb{R}$ such that
$b(\cdot,x_0)\in L^{\kappa_1}_+([0,T];\mathbb{R}), \sigma(\cdot,y_0)\in L^{2\kappa_2}_+([0,T];\mathbb{R})$.
Denote by $X^x$ and $X^y$ the solution of \eqref{2.3} with initial condition $x$ and $y$ respectively. Then, for any $p\in A_H$, there exists constant
$L_p>0$ such that
\beqnn
\sup\limits_{t\in[0,T]}\mathbb{E}|X_t^x-X_t^y|^p\leq L_p|x-y|^p.
\eeqnn
\elemma

\bremark\label{r4.2}
If we consider the case $p=2$, Lemma \ref{l4.2} reduces to \cite[Theorem 3.3]{COutin&Decreu}.
\eremark

Next we will concern the existence of the derivative process w.r.t. the initial data.
\bproposition\label{p4.1}
Suppose that $b$ and $\sigma$ are both differentiable w.r.t. their second variables such that
 $\partial b(t,\cdot)$ and $\partial\sigma(t,\cdot)$ are bounded and uniformly continuous uniformly w.r.t. their first variable $t$.
 Then, for each $y\in\mathbb{R}$, $(\nabla_yX_t^x)_{0\leq t\leq T}$ exists and is the unique solution to the equation
\beqnn
Y_t=y+\int_0^tK_H(t,s)\partial b(s,\cdot)(X_s^x)Y_s{\rm d}s+\int_0^tK_H(t,s)\partial\sigma(s,\cdot)(X_s^x)Y_s{\rm d}W_s,
\eeqnn
where $X_\cdot^x$ is the solution of \eqref{2.3}.
\eproposition

\emph{Proof.} Using the Picard iteration argument introduced in \cite[ Theorem 3.1]{COutin&Decreu}, we can easily show that the above equation
has a unique solution $(Y_t)_{t\in[0,T]}$ and moreover, $\sup_{t\in[0,T]}\mathbb{E}|Y_t|^p<\infty$ holds for any $p\in A_H$. \\
For $\epsilon >0$, let $Z_t^\epsilon=X_t^{x+\epsilon y}-X_t^x-\epsilon Y_t,\ t\in[0,T]$. To complete the proof, it suffices to prove
\beqnn
\lim\limits_{\epsilon\rightarrow 0}\mathbb{E}\frac{|Z_t^\epsilon|^2}{\epsilon^2}=0,\ \forall t\in[0,T].
\eeqnn
To this end, we see that, for any $t\in[0,T]$,
\beqnn
Z_t^\epsilon\ar=\ar\int_0^tK_H(t,s)\left(b(s,X_s^{x+\epsilon y})-b(s,X_s^x)-\epsilon\partial b(s,\cdot)(X_s^x)Y_s\right){\rm d}s\cr
\ar+\ar\int_0^tK_H(t,s)\left(\sigma(s,X_s^{x+\epsilon y})-\sigma(s,X_s^x)-\epsilon\partial \sigma(s,\cdot)(X_s^x)Y_s\right){\rm d}W_s.
\eeqnn
Therefore, by the H\"{o}lder inequality and the Burkholder-Davis-Gundy inequality, there is some constant $C_2$ such that
\beqlb\label{4.1}
\mathbb{E}|Z_t^\epsilon|^2\ar\leq\ar2T\int_0^tK_H^2(t,s)\left|b(s,X_s^{x+\epsilon y})-b(s,X_s^x)-\epsilon\partial b(s,\cdot)(X_s^x)Y_s\right|^2{\rm d}s\cr
\ar+\ar2C_2\int_0^tK_H^2(t,s)\left|\sigma(s,X_s^{x+\epsilon y})-\sigma(s,X_s^x)-\epsilon\partial \sigma(s,\cdot)(X_s^x)Y_s\right|^2{\rm d}s\cr
\ar=:\ar2T\int_0^tK_H^2(t,s)J_1(s)^2{\rm d}s+2C_2\int_0^tK_H^2(t,s)J_2(s)^2{\rm d}s.
\eeqlb
Next we are to estimate $J_1(s)$ and $J_2(s)$. Let us define, for each $\delta\geq0$,
\beqnn
\alpha(\delta)=\sup\limits_{|x-y|\leq\delta}\sup\limits_{s\in[0,T]}\left(|\partial b(s,\cdot)(x)-\partial b(s,\cdot)(y)|
+|\partial \sigma(s,\cdot)(x)-\partial \sigma(s,\cdot)(y)|\right).
\eeqnn
It is clear from the assumptions on the coefficients $b$ and $\sigma$ that $\alpha(\infty)<\infty$ and $\alpha(\delta)\downarrow0$ as $\delta\downarrow0$.
As a consequence, we derive that,
\beqnn
\delta^2\alpha^2(\delta)=\delta^2\alpha^2(\delta){\rm I}_{\{\delta\leq\sqrt{\epsilon}\}}+\delta^2\alpha^2(\delta){\rm I}_{\{\delta>\sqrt{\epsilon}\}}
\leq\delta^2\alpha^2(\sqrt{\epsilon})+\frac{\delta^q\alpha^2(\infty)}{\epsilon^{\frac{q-2}{2}}},
\eeqnn
where $q$ is chosen such that $2<q<\frac{1}{H_0}$.\\
Note that, by the mean value theorem, we get
\beqnn
J_1(s)=|\left(\partial b(s,\cdot)(\zeta_1)-\partial b(s,\cdot)(X_s^x)\right)(X_s^{x+\epsilon y}-X_s^x)
    +\partial b(s,\cdot)(X_s^x)Z_s^\epsilon|
\eeqnn
and
\beqnn
J_2(s)=|\left(\partial\sigma(s,\cdot)(\zeta_2)-\partial\sigma(s,\cdot)(X_s^x)\right)(X_s^{x+\epsilon y}-X_s^x)
    +\partial\sigma(s,\cdot)(X_s^x)Z_s^\epsilon|,
\eeqnn
where $\zeta_i=\theta_iX_s^x+(1-\theta_i)X_s^{x+\epsilon y}, \theta_i\in(0,1), i=1,2$.\\
Hence, we conclude that
\beqlb\label{4.2}
J_1(s)^2+J_2(s)^2\ar\leq\ar(J_1(s)+J_2(s))^2\leq2\alpha^2(|X_s^{x+\epsilon y}-X_s^x|)|X_s^{x+\epsilon y}-X_s^x|^2+2M|Z_s^\epsilon|^2\cr
\ar\leq\ar2\alpha^2(\sqrt{\epsilon})|X_s^{x+\epsilon y}-X_s^x|^2+\frac{\alpha^2(\infty)}{\epsilon^{\frac{q-2}{2}}}|X_s^{x+\epsilon y}-X_s^x|^q
+2M|Z_s^\epsilon|^2,
\eeqlb
where $M=\sup_{s\in[0,T]}(|\partial b(s,\cdot)|+|\partial\sigma(s,\cdot)|)^2$.\\
Now we turn to the estimate of $\mathbb{E}|Z_t^\epsilon|^2$. Substituting \eqref{4.2} into \eqref{4.1} and noting $q\in A_H$, we have, by Lemma \ref{l4.2},
\beqnn
\mathbb{E}|Z_t^\epsilon|^2\ar\leq\ar2(T\vee C_2)T^{2H}\left(2L_2|y|^2\alpha^2(\sqrt{\epsilon})\epsilon^2+L_q|y|^q\alpha^2(\infty)\epsilon^{\frac{q}{2}+1}\right)
+4(T\vee C_2)M\int_0^tK_H^2(t,s)\mathbb{E}|Z_s^\epsilon|^2{\rm d}s\cr
\ar=\ar:C(\epsilon)+C_3\int_0^tK_H^2(t,s)\mathbb{E}|Z_s^\epsilon|^2{\rm d}s.
\eeqnn
We set $K_1^2(t,s)=K_H^2(t,s), K_{n+1}^2(t,s)=\int_s^tK_1^2(t,r)K_n^2(r,s){\rm d}r, \forall s,t\in[0,T], n\geq1$, and identify the operator $K_n^2$ with its kernel,
$K_0^2:\equiv1$. \\
Then by induction, we deduce that
\beqnn
\mathbb{E}|Z_t^\epsilon|^2\ar\leq\ar C(\epsilon)(1+C_3(K_1^21)(t))+C_3^2\int_0^t\int_0^sK_H^2(t,s)K_H^2(s,r)\mathbb{E}|Z_r^\epsilon|^2{\rm d}r{\rm d}s\cr
\ar\leq\ar\cdot\cdot\cdot\leq C(\epsilon)\sum_{i=0}^nC_3^i(K_i^21)(t)+C_3^{n+1}\sup\limits_{u\in[0,T]}\mathbb{E}|Z_u^\epsilon|^2(K_{n+1}^21)(t).
\eeqnn
Recall that \cite[Lemma 3.3]{COutin&Decreu} states that $\sum_{i=0}^\infty\sup_{0\leq t\leq T}(K_i^21)(t)z^i<\infty,\ \forall z\in\mathbb{C}$.
Therefore, letting $n\rightarrow\infty$, it follows that
\beqnn
\mathbb{E}|Z_t^\epsilon|^2\ar\leq\ar C(\epsilon)\sum_{i=0}^\infty C_3^i\sup_{t\in[0,T]}(K_i^21)(t).
\eeqnn
Observing that
$\lim\limits_{\epsilon\rightarrow0}\frac{C(\epsilon)}{\epsilon^2}=0$,
the proof is finished.

Now we are in position to prove Theorem \ref{4.1}.\\ \\
\emph{Proof of Theorem \ref{4.1}.}
Note that if there exists $\xi\in$ Dom$\delta_H$ such that
\beqlb\label{4.3}
\langle D_HX_T^x, \xi\rangle_\mathcal {H}=\nabla_yX_T^x,\  a.s.,
\eeqlb
then for each $f\in C_b^1(\mathbb{R})$,
\beqnn
\nabla_yP_Tf(x)\ar=\ar\nabla_y\mathbb{E}f(X_T^x)=\mathbb{E}\nabla_yf(X_T^x)=\mathbb{E}(f'(X_T^x)\nabla_yX_T^x)\cr
\ar=\ar\mathbb{E}(f'(X_T^x)\langle D_HX_T^x, \xi\rangle_\mathcal {H})=\mathbb{E}(\langle D_Hf(X_T^x), \xi\rangle_\mathcal {H}).
\eeqnn
Applying the integration by parts formula for $D_H$, i.e. the definition of $\delta_H$, we get
\beqnn
\nabla_yP_Tf(x)=\mathbb{E}(f(X_T^x)\delta_H\xi)=\mathbb{E}\left(f(X_T^x)\int_0^T\dot{\xi}_s\delta_HB_s^H\right).
\eeqnn
Furthermore, if $\dot{\xi}\in L_a^2([0,T]\times\Omega)$, then $\nabla_yP_Tf(x)=\mathbb{E}\left(f(X_T^x)\int_0^T\dot{\xi}_s{\rm d}W_s\right)$.\\
Based on the analysis above, we know that, to complete the proof, it suffices to find a $\xi=K_H\dot{\xi}$ such that $\dot{\xi}\in L_a^2([0,T]\times\Omega)$
and \eqref{4.3} holds.\\
Let
\beqnn
 \dot{\xi}_s=\sigma(s)^{-1}\left(\left(1+\int_0^sK_H(s,r)u'(r){\rm d}r\right)\partial b(s,\cdot)(X_s^x)-u'(s)\right)y,
\eeqnn
where $u\in C^1([0,T];\mathbb{R})$ such that $1+\int_0^TK_H(T,r)u'(r){\rm d}r=0$.
Obviously, $\dot{\xi}$ constructed above is in $L_a^2([0,T]\times\Omega)$.
Next consider the following equation
\beqlb\label{4.4}
Z_t=y+\int_0^tK_H(t,s)\partial b(s,\cdot)(X_s^x)Z_s{\rm d}s-\int_0^tK_H(t,s)\sigma(s)\dot{\xi}_s{\rm d}s.
\eeqlb
By the assumption, it is clear that \eqref{4.4} has a unique solution $Z$. On one hand, observe that $Y_t:=(1+\int_0^tK_H(t,r)u'(r){\rm d}r)y$ solves \eqref{4.4}.
On the other hand, since $\partial\sigma(s,\cdot)=0, \forall s\in[0,T]$, Lemma \ref{l4.1} together with Proposition \ref{p4.1} implies that
$(\nabla_yX_t^x-\langle D_HX_t^x, \xi\rangle_\mathcal {H})_{t\in[0,T]}$ is also a solution of \eqref{4.4}. As a consequence,
$Y_t=\nabla_yX_t^x-\langle D_HX_t^x, \xi\rangle_\mathcal {H}, \forall t\in[0,T]$ holds. Due to $Y_T=0$, it follows that
$\langle D_HX_T^x, \xi\rangle_\mathcal {H}=\nabla_yX_T^x$. Therefore, the proof is complete.

\bremark\label{r4.3}
If we take $u'(t)=-\frac{C_H}{T}t^{\frac{1}{2}-H}$, then
by the proof of Theorem \ref{t3.1}, we know that $1+\int_0^TK_H(T,r)u'(r){\rm d}r=0$ holds and
the result of Theorem \ref{t4.1} can be expressed as
\beqnn
\nabla_y P_T f(x)=\mathbb{E}\left[f(X_T^x)\int_0^T\sigma(s)^{-1}\left((T-s)\partial b(s,\cdot)(X_s^x)+C_Hs^{\frac{1}{2}-H}\right)\frac{y}{T}{\rm d}W_s\right].
\eeqnn
In particular, when $H=\frac{1}{2}$, we obtain a version of relation above that is an extension of \cite[Theorem 3.1]{Wang12a},
in which the coupling argument is used.
\eremark

Next we will state some applications of the derivative formula obtained above. More precisely, explicit gradient estimate, Harnack inequality and log-Harnack
inequality are presented. That is

\bcorollary\label{c4.1}
Assume that (H2) holds and set $C(T,K_3,K_4,H)=2K_4^2\left(\frac{K_3^2T}{3}+\frac{C_H^2}{(2-2H)T^{2H}}\right)$.
\begin{itemize}
\item[{\rm (1)}]
 For any  $f\in\mathscr{B}_b(\mathbb{R})$, we get
\beqnn
|\nabla_yP_Tf(x)|^2\leq C(T,K_3,K_4,H)|y|^2P_Tf^2(x),
\eeqnn
i.e., $|\nabla_yP_Tf(x)|$ is bounded above by $f$.
  Moreover, for all $\delta>0$ and positive $f\in\mathscr{B}_b(\mathbb{R})$,
\beqlb\label{4.5}
 |\nabla_yP_Tf(x)|\leq\delta\left[P_T(f\log f)-(P_Tf)(\log P_Tf)\right](x)+\frac{C(T,K_3,K_4,H)}{\delta}|y|^2P_Tf(x).
\eeqlb
 \par
 \item[{\rm (2)}]
For any non-negative $f\in\mathscr{B}_b(\mathbb{R})$ and $p>1$, the following Harnack inequality holds:
\beqlb\label{4.6}
(P_Tf(x))^p\leq P_Tf^p(y){\rm exp}\left[\frac{p}{p-1}C(T,K_3,K_4,H)|x-y|^2\right],\ x, y\in\mathbb{R}.
\eeqlb
 As a consequence, the log-Harnack inequality
\beqlb\label{4.7}
P_T(\log f)(x)\leq\log P_Tf(y)+C(T,K_3,K_4,H)|x-y|^2
\eeqlb
holds for any positive $f\in\mathscr{B}_b(\mathbb{R})$,
and
$P_T$ is strong Feller, i.e. for each $x\in\mathbb{R}$,
\beqnn
\lim\limits_{y\rightarrow x}P_Tf(y)=P_Tf(x).
\eeqnn
\item[{\rm (3)}]
Let $\mu$ be $P_T$ sub-invariant, i.e., $\mu$ is a probability measure on $\mathbb{R}$ such that
$\int_\mathbb{R}P_Tf{\rm d}\mu\leq\int_\mathbb{R}f{\rm d}\mu$ for all $f\in\mathscr{B}_b(\mathbb{R}), f\geq0$. Then the entropy-cost inequality
\beqlb\label{4.8}
\mu(P_T^\ast f\log P_T^\ast f)\leq C(T,K_3,K_4,H)W_2^d(f\mu,\mu)^2, f\geq0, \mu(f)=1,
\eeqlb
holds for the adjoint operator $P_T^\ast$ of $ P_T$ in $L^2(\mu)$, where $d(x,y)=|x-y|$.
\end{itemize}
\ecorollary
\emph{Proof.}
Let $u'(t)=-\frac{C_H}{T}t^{\frac{1}{2}-H}$ and define
$M_T=\int_0^T\sigma(s)^{-1}\left((T-s)\partial b(s,\cdot)(X_s^x)+C_Hs^{\frac{1}{2}-H}\right)\frac{y}{T}{\rm d}W_s$.
By the hypotheses on the coefficients, we derive that
\beqnn
\langle M\rangle_T=\int_0^T\left|\sigma(s)^{-1}\left((T-s)\partial b(s,\cdot)(X_s^x)+C_Hs^{\frac{1}{2}-H}\right)\frac{y}{T}\right|^2{\rm d}s
\leq C(T,K_3,K_4,H)|y|^2,
\eeqnn
where $C(T,K_3,K_4,H)=2K_4^2\left(\frac{K_3^2T}{3}+\frac{C_H^2}{(2-2H)T^{2H}}\right)$.\\
Hence, it follows from the H\"{o}lder inequality that
\beqnn
|\nabla_yP_Tf(x)|^2\leq\mathbb{E}\langle M\rangle_TP_Tf^2(x)\leq C(T,K_3,K_4,H)|y|^2P_Tf^2(x).
\eeqnn
Combining the derivative formula with the Young inequality yield that, for any positive $f\in\mathscr{B}_b(\mathbb{R})$ and $\delta>0$,
\beqlb\label{4.9}
|\nabla_yP_Tf(x)|\leq\delta\left[P_T(f\log f)-(P_Tf)(\log P_Tf)\right]+\delta\log\mathbb{E}\exp\left[\frac{M_T}{\delta}\right]P_Tf.
\eeqlb
 Observe that
\beqlb\label{4.10}
\mathbb{E}\exp\left[\frac{M_T}{\delta}\right]\leq\left(\mathbb{E}\exp\left[\frac{2\langle M\rangle_T}{\delta^2}\right]\right)^\frac{1}{2}
\leq\exp\left[\frac{C(T,K_3,K_4,H)}{\delta^2}|y|^2\right].
\eeqlb
Substituting \eqref{4.10} into \eqref{4.9} implies \eqref{4.5}.
In the sprit of \cite[Corollary 1.2]{Wang13}, \eqref{4.6} follows from \eqref{4.5}. Since $\mathbb{R}$ is a length space, then according to
\cite[Proposition 2.2]{Wang10}, \eqref{4.6} implies \eqref{4.7}. The strong Feller property follows from  \eqref{4.6}, due to the same proof of
\cite[Proposition 4.1]{Da}. Finally, \eqref{4.8} can be proved as the proof of \cite[Corollary 1.2]{Rockner} or \cite[Corollary 3.6]{Fan}.

\bremark\label{r4.4}
Making use of the Harnack type inequalities, one can compare the values of a reference function at different points, while in the shift Harnack type inequalities presented in Theorem \ref{3.1}, instead of initial points, a reference function is shifted.
Besides,  the (resp. shift) Harnack type inequalities allow us to compare the measure $P_T(x,\cdot)$ with some invariant probability
measure associated with a certain semigroup (resp. the Lebesgue measure), where $P_T(x,\cdot)$ is the transition probability for $P_T$.
One can see \cite{Wang12b} for more applications of  the shift Harnack type inequalities.
\eremark

\section{Transportation inequalities}

\setcounter{equation}{0}
In this section we will discuss the Talagrand type transportation cost inequalities for the law of the solution of \eqref{2.3} w.r.t. the uniform
distance $d_\infty$ and the $L^2$-distance $d_2$ on the path space $C([0,T];\mathbb{R})$.
To the end, we introduce the following assumption: (H3)
\begin{itemize}
 \item[(i)] there exists constant $K_6(>0)$ such that
 \beqnn
 |b(t,x)-b(t,y)|+|\sigma(t,x)-\sigma(t,y)|\leq K_6|x-y|,\ \forall t\in[0,T],\ x,y\in\mathbb{R};
 \eeqnn
\par
 \item[(ii)] $\|\sigma\|_\infty:=\sup\limits_{0\leq t\leq T}\sup\limits_{x\in\mathbb{R}}|\sigma(t,x)|<\infty, b(\cdot,0)\in L^{\kappa_1}_+([0,T];\mathbb{R})$.
\end{itemize}

Let us start by proving the following proposition which is crucial for the proof of Theorem \ref{t5.1} below.
\bproposition\label{p5.1}
Let $H>\frac{1}{2}$ and $\tau$ be an $(\mathscr{F}_t)$-stopping time. Assume that $\phi$ is an adapted stochastic process satisfying $\mathbb{E}\int_0^T|\phi_t|^p{\rm d}t<\infty$ for some $p\geq2$. Then,
there holds the maximal inequality
\beqnn
\mathbb{E}\left(\sup\limits_{0\leq t\leq{T\wedge\tau}}\left|\int_0^tK_H(t,s)\phi(s){\rm d}W_s\right|^p\right)\leq C(p)\mathbb{E}\int_0^{T\wedge\tau}|\phi_t|^p{\rm d}t,
\eeqnn
where $C(p)$ is a positive constant depending on $p$.
\eproposition
\emph{Proof.}
Recall that, for $H\in(0,1), K_H(t,s)$ is the kernel
\beqnn
K_H(t,s)=\alpha_H(t-s)^{H-\frac{1}{2}}+\alpha_H\left(\frac{1}{2}-H\right)\int_s^t(r-s)^{H-\frac{3}{2}}\left(1-\left(\frac{s}{r}\right)^{\frac{1}{2}-H}\right){\rm d}r,
\eeqnn
where $\alpha_H=\left(\frac{2H\Gamma(3/2-H)}{\Gamma(H+1/2)\Gamma(2-2H)}\right)^{1/2}$.\\
From this relation, we get
\beqnn
\frac{\partial K_H(t,s)}{\partial t}=\alpha_H\left(H-\frac{1}{2}\right)\left(\frac{s}{t}\right)^{\frac{1}{2}-H}(t-s)^{H-\frac{3}{2}}.
\eeqnn
When $H>\frac{1}{2}$, the kernel $K_H(t,s)$ can reformulate as (for instance, see \cite{Alos} and references therein)
\beqnn
K_H(t,s)=\alpha_H\left(H-\frac{1}{2}\right)s^{\frac{1}{2}-H}\int_s^tr^{H-\frac{1}{2}}(r-s)^{H-\frac{3}{2}}{\rm d}r
=:\bar{\alpha}_Hs^{\frac{1}{2}-H}\int_s^tr^{H-\frac{1}{2}}(r-s)^{H-\frac{3}{2}}{\rm d}r.
\eeqnn
Therefore, we have
\beqnn
\int_0^tK_H(t,s)\phi(s){\rm d}W_s=\bar{\alpha}_H\int_0^ts^{\frac{1}{2}-H}\int_s^tr^{H-\frac{1}{2}}(r-s)^{H-\frac{3}{2}}{\rm d}r\phi(s){\rm d}W_s.
\eeqnn
To exchange the integration of the right-hand side of the above expression, one need to show that the integrand fulfills the conditions of the stochastic Fubini theorem (see \cite[Theorem 4.18]{Da92}).
Actually, choosing $\epsilon\in(0,\frac{1}{2})$ such that $H>\frac{1+\epsilon}{2}$ and using the H\"{o}lder inequality and the Young inequality, we obtain
\beqnn
\ar\ar\int_0^tr^{H-\frac{1}{2}}\left(\mathbb{E}\int_0^rs^{1-2H}(r-s)^{2H-3}\phi^2(s){\rm d}s\right)^\frac{1}{2}{\rm d}r\cr
\ar\leq\ar\left(\int_0^tr^{-1+2\epsilon}{\rm d}r\right)^{\frac{1}{2}}\left(\int_0^tr^{2(H-\epsilon)}\mathbb{E}\int_0^rs^{1-2H}(r-s)^{2H-3}\phi^2(s){\rm d}s{\rm d}r\right)^{\frac{1}{2}}\cr
\ar\leq\ar\left(\frac{T^{2\epsilon}}{2\epsilon}\right)^\frac{1}{2}\left(\int_0^tr^{4H-2\epsilon-3}{\rm d}r
\mathbb{E}\int_0^tr^{1-2\epsilon}\phi^2(r){\rm d}r\right)^{\frac{1}{2}}\cr
\ar\leq\ar\left(\frac{T^{4H-2\epsilon-1}}{4\epsilon(2H-\epsilon-1)}\right)^\frac{1}{2}\left(\mathbb{E}\int_0^T\phi^2(r){\rm d}r\right)^{\frac{1}{2}},
\eeqnn
which is finite due to hypothesis on $\phi$.\\
So, the stochastic Fubini theorem implies
\beqnn
\int_0^tK_H(t,s)\phi(s){\rm d}W_s=\bar{\alpha}_H\int_0^tr^{H-\frac{1}{2}}\int_0^rs^{\frac{1}{2}-H}(r-s)^{H-\frac{3}{2}}\phi(s){\rm d}W_s{\rm d}r.
\eeqnn
Taking $\theta\in(0,\frac{1}{2})$ such that $H>\frac{1+\theta}{2}$ and applying the H\"{o}lder inequality and the Young inequality, we obtain
\beqnn
\ar\ar\mathbb{E}\left(\sup\limits_{0\leq t\leq {T\wedge\tau}}\left|\int_0^tr^{H-\frac{1}{2}}\int_0^rs^{\frac{1}{2}-H}(r-s)^{H-\frac{3}{2}}\phi(s){\rm d}W_s{\rm d}r\right|^p\right)\cr
\ar\leq\ar\left(\int_0^Tr^{(\frac{1}{p}-1+\theta)\frac{p}{p-1}}{\rm d}r\right)^{p-1}
\int_0^Tr^{(H+\frac{1}{2}-\frac{1}{p}-\theta)p}\mathbb{E}\left(\left|\int_0^rs^{\frac{1}{2}-H}(r-s)^{H-\frac{3}{2}}\phi(s){\rm d}W_s\right|^p\cdot{\rm I}_{[0,T\wedge\tau]}(r)\right){\rm d}r\cr
\ar\leq\ar\left(\int_0^Tr^{(\frac{1}{p}-1+\theta)\frac{p}{p-1}}{\rm d}r\right)^{p-1}
\int_0^Tr^{(H+\frac{1}{2}-\frac{1}{p}-\theta)p}\mathbb{E}\left|\int_0^{r\wedge\tau}s^{\frac{1}{2}-H}(r-s)^{H-\frac{3}{2}}\phi(s){\rm d}W_s\right|^p{\rm d}r\cr
\ar\leq\ar\left(\frac{p-1}{\theta p}\right)^{p-1}T^{\theta p}\int_0^Tr^{(H+\frac{1}{2}-\frac{1}{p}-\theta)p}\mathbb{E}\left(\int_0^{r\wedge\tau}s^{1-2H}(r-s)^{2H-3}\phi(s)^2{\rm d}s\right)^\frac{p}{2}{\rm d}r\cr
\ar=\ar\left(\frac{p-1}{\theta p}\right)^{p-1}T^{\theta p}\int_0^Tr^{(H+\frac{1}{2}-\frac{1}{p}-\theta)p}\mathbb{E}\left(\int_0^Ts^{1-2H}(r-s)^{2H-3}\phi(s)^2{\rm I}_{\{s\leq r\}}{\rm I}_{\{s\leq T\wedge\tau\}}{\rm d}s\right)^\frac{p}{2}{\rm d}r\cr
\ar\leq\ar\left(\frac{p-1}{\theta p}\right)^{p-1}T^{\theta p}\left(\int_0^Tr^{(H+\frac{1}{2}-\frac{1}{p}-\theta)p+2H-3}{\rm d}r\right)^\frac{p}{2}
\mathbb{E}\int_0^Tr^{(H+\frac{1}{2}-\frac{1}{p}-\theta)p+(1-2H)\frac{p}{2}}|\phi_r|^p{\rm I}_{\{r\leq{T\wedge\tau}\}}{\rm d}r\cr
\ar\leq\ar\left(\frac{p-1}{\theta p}\right)^{p-1}\left(\frac{1}{(H+\frac{1}{2}-\theta)p+2H-3}\right)^\frac{p}{2}T^{(H+\frac{1}{2}-\theta)\frac{p^2}{2}+(H-\frac{1}{2})p-1}
\mathbb{E}\int_0^{T\wedge\tau}|\phi_r|^p{\rm d}r,
\eeqnn
which yields the desired result.

We now prove the following main result in this section.
\btheorem\label{t5.1}
Let $H>\frac{1}{2}$. Assume (H3) and let $\mathbb{P}_x$ be the law of the solution of \eqref{2.3} with the initial point $x$ on the path space $C([0,T];\mathbb{R})$. Then,
$\mathbb{P}_x$ satisfies the transportation cost inequalities on the metric space $C([0,T];\mathbb{R})$. More precisely,
\begin{itemize}
 \item[(1)] $\mathbb{P}_x\in T_2(\alpha(T,H)|d_\infty)$, where $\alpha(T,H)=3(\|\sigma\|_\infty T^H)^2{\rm e}^{3K_6^2T(T^{2H}+C(2))}$,
\par
 \item[(2)] $\mathbb{P}_x\in T_2(\beta(T,H)|d_2)$, where $\beta(T,H)=3(\|\sigma\|_\infty T^H)^2\frac{{\rm e}^{3K_6^2T(T^{2H}+C(2))}-1}{3K_6^2(T^{2H}+C(2))}$.
\end{itemize}
\etheorem
\emph{Proof.}
Let $\mathbb{Q}$ be a probability measure on $C([0,T];\mathbb{R})$ such that $\mathbb{Q}\ll\mathbb{P}_x$. Clearly, to  prove the desired result, we only need
to consider the case $\mathbb{H}(\mathbb{Q}|\mathbb{P}_x)<\infty$. The proof will divide into two steps.\\
Step 1. The part follows the arguments of \cite{Djel}. Let $\bar{\mathbb{Q}}=\frac{{\rm d}\mathbb{Q}}{{\rm d}\mathbb{P}_x}(X_\cdot)\mathbb{P}$. Note that
\beqnn
\int_\Omega\frac{{\rm d}\mathbb{Q}}{{\rm d}\mathbb{P}_x}(X_\cdot){\rm d}\mathbb{P}
=\int_{C([0,T];\mathbb{R})}\frac{{\rm d}\mathbb{Q}}{{\rm d}\mathbb{P}_x}(\gamma){\rm d}\mathbb{P}_x(\gamma)=\mathbb{Q}(C([0,T];\mathbb{R}))=1,
\eeqnn
and
\beqnn
\int_{C([0,T];\mathbb{R})}\frac{{\rm d}\mathbb{Q}}{{\rm d}\mathbb{P}_x}\log\left(\frac{{\rm d}\mathbb{Q}}{{\rm d}\mathbb{P}_x}\right){{\rm d}\mathbb{P}_x}
=\int_\Omega\frac{{\rm d}\mathbb{Q}}{{\rm d}\mathbb{P}_x}(X^x_\cdot)\log\left(\frac{{\rm d}\mathbb{Q}}{{\rm d}\mathbb{P}_x}(X^x_\cdot)\right){{\rm d}\mathbb{P}}
=\int_\Omega\frac{{\rm d}\bar{\mathbb{Q}}}{{\rm d}\mathbb{P}}\log\left(\frac{{\rm d}\bar{\mathbb{Q}}}{{\rm d}\mathbb{P}}\right){{\rm d}\mathbb{P}},
\eeqnn
that is, $\bar{\mathbb{Q}}$ is a probability measure on $(\Omega,\mathcal {F})$ and $\mathbb{H}(\mathbb{Q}|\mathbb{P}_x)=\mathbb{H}(\bar{\mathbb{Q}}|\mathbb{P})$.\\
According to the proof of \cite[Theorem 5.6]{Djel}, there is a predictable process $(u_t)_{0\leq t\leq T}$ such that
\beqnn
\mathbb{H}(\mathbb{Q}|\mathbb{P}_x)=\mathbb{H}(\bar{\mathbb{Q}}|\mathbb{P})=\frac{1}{2}\mathbb{E}_{\bar{\mathbb{Q}}}\int_0^T|u_t|^2{\rm d}t
\eeqnn
and the process
\beqnn
\bar{W}_t:=W_t-\int_0^tu_s{\rm d}s
\eeqnn
is a Brownian motion under $\bar{\mathbb{Q}}$, where $\mathbb{E}_{\bar{\mathbb{Q}}}$ is the expectation taken for the probability measure $\bar{\mathbb{Q}}$.
As a consequence, the process $(\bar{B}^H_t)_{0\leq t\leq T}$ defined by
\beqnn
\bar{B}^H_t=\int_0^tK_H(t,s){\rm d}\bar{W}_s=B^H_t-(K_Hu)(t)
\eeqnn
is a $\bar{\mathbb{Q}}$-fractional Brownian motion associated with $\bar{W}$.\\
Step 2.
From step 1, we can reformulate \eqref{2.3} as
\beqlb\label{5.1}
X_t=x+\int_0^tK_H(t,s)\left(b(s,X_s)+\sigma(s,X_s)u_s\right){\rm d}s+\int_0^tK_H(t,s)\sigma(s,X_s){\rm d}\bar{W}_s.
\eeqlb
Noting that, for any bounded measurable function $F$ on $C([0,T];\mathbb{R})$,
\beqnn
\mathbb{E}_{\bar{\mathbb{Q}}}(F(X_\cdot))=\mathbb{E}\left(\frac{{\rm d}\mathbb{Q}}{{\rm d}\mathbb{P}_x}(X_\cdot)F(X_\cdot)\right)
=\int_{C([0,T];\mathbb{R})}\frac{{\rm d}\mathbb{Q}}{{\rm d}\mathbb{P}_x}(\gamma)F(\gamma){{\rm d}\mathbb{P}_x}(\gamma)
=\mathbb{Q}(F),
\eeqnn
it follows that the law of $X_\cdot$ under $\bar{\mathbb{Q}}$ is $\mathbb{Q}$. On the other hand, we consider the following equation
\beqlb\label{5.2}
Y_t=x+\int_0^tK_H(t,s)b(s,Y_s){\rm d}s+\int_0^tK_H(t,s)\sigma(s,Y_s){\rm d}\bar{W}_s.
\eeqlb
As $\bar{W}$ is the Brownian motion under $\bar{\mathbb{Q}}$, we easily know that the law of $Y_\cdot$ under
$\bar{\mathbb{Q}}$ is $\mathbb{P}_x$.
Therefore, the law of $(X,Y)$ under $\bar{\mathbb{Q}}$ is a coupling of $(\mathbb{Q},\mathbb{P}_x)$ and moreover, we get
\beqnn
\ar\ar W_2^{d_\infty}(\mathbb{Q},\mathbb{P}_x)^2\leq\mathbb{E}_{\bar{\mathbb{Q}}}d_\infty(X,Y)^2=\mathbb{E}_{\bar{\mathbb{Q}}}\left(\sup\limits_{0\leq t\leq T}|X_t-Y_t|^2\right),\cr
\ar\ar W_2^{d_2}(\mathbb{Q},\mathbb{P}_x)^2\leq\mathbb{E}_{\bar{\mathbb{Q}}}d_2(X,Y)^2=\mathbb{E}_{\bar{\mathbb{Q}}}\left(\int_0^T|X_t-Y_t|^2{\rm d}t\right).
\eeqnn
Combining \eqref{5.1} with \eqref{5.2}, we have
\beqnn
X_t-Y_t\ar=\ar\int_0^tK_H(t,s)\left(b(s,X_s)-b(s,Y_s)\right){\rm d}s+\int_0^tK_H(t,s)\sigma(s,X_s)u_s{\rm d}s\cr
\ar+\ar\int_0^tK_H(t,s)\left(\sigma(s,X_s)-\sigma(s,Y_s)\right){\rm d}\bar{W}_s.
\eeqnn
Now, for $n\in\mathbb{N}$, define the stopping time
\beqnn
\tau_n:=\inf\{t>0,|X_t-Y_t|\geq n\}.
\eeqnn
Obviously, $\tau_n\uparrow\infty$ as $n$ goes to $\infty$.
Applying Proposition \ref{p5.1} and the H\"{o}lder inequality, we derive that
\beqnn
\mathbb{E}_{\bar{\mathbb{Q}}}\left(\sup\limits_{0\leq t\leq T}\left|(X_t-Y_t){\rm I}_{\{t\leq\tau_n\}}\right|^2\right)
\ar\leq\ar3\mathbb{E}_{\bar{\mathbb{Q}}}\left(\sup\limits_{0\leq t\leq T}\int_0^{t\wedge\tau_n}\left|K_H(t,s)(b(s,X_s)-b(s,Y_s))\right|{\rm d}s\right)^2\cr
\ar+\ar3\mathbb{E}_{\bar{\mathbb{Q}}}\left(\sup\limits_{0\leq t\leq T}\int_0^{t\wedge\tau_n}\left|K_H(t,s)\sigma(s,X_s)u_s\right|{\rm d}s\right)^2\cr
\ar+\ar3\mathbb{E}_{\bar{\mathbb{Q}}}\left(\sup\limits_{0\leq t\leq T}
\left|{\rm I}_{\{t\leq\tau_n\}}\cdot\int_0^tK_H(t,s)\left(\sigma(s,X_s)-\sigma(s,Y_s)\right){\rm d}\bar{W}_s\right|^2\right)\cr
\ar\leq\ar3K_6^2\mathbb{E}_{\bar{\mathbb{Q}}}\left(\sup\limits_{0\leq t\leq T}\int_0^{t\wedge\tau_n}|K_H(t,s)|\cdot|X_s-Y_s|{\rm d}s\right)^2\cr
\ar+\ar3\|\sigma\|_\infty^2\mathbb{E}_{\bar{\mathbb{Q}}}\left(\sup\limits_{0\leq t\leq T}\int_0^{t\wedge\tau_n}\left|K_H(t,s)u_s\right|{\rm d}s\right)^2\cr
\ar+\ar3C(2)\mathbb{E}_{\bar{\mathbb{Q}}}\left(\int_0^{T\wedge\tau_n}\left|\sigma(s,X_s)-\sigma(s,Y_s)\right|^2{\rm d}s\right)\cr
\ar\leq\ar3(\|\sigma\|_\infty T^H)^2\mathbb{E}_{\bar{\mathbb{Q}}}\int_0^Tu_s^2{\rm d}s\cr
\ar+\ar3K_6^2(T^{2H}+C(2))\mathbb{E}_{\bar{\mathbb{Q}}}\left(\int_0^{T\wedge\tau_n}\left|X_s-Y_s\right|^2{\rm d}s\right)\cr
\ar=\ar3(\|\sigma\|_\infty T^H)^2\mathbb{E}_{\bar{\mathbb{Q}}}\int_0^Tu_s^2{\rm d}s\cr
\ar+\ar3K_6^2(T^{2H}+C(2))\mathbb{E}_{\bar{\mathbb{Q}}}\left(\int_0^T\left|(X_s-Y_s){\rm I}_{\{s\leq\tau_n\}}\right|^2{\rm d}s\right)\cr
\ar\leq\ar3(\|\sigma\|_\infty T^H)^2\mathbb{E}_{\bar{\mathbb{Q}}}\int_0^Tu_s^2{\rm d}s\cr
\ar+\ar3K_6^2(T^{2H}+C(2))\int_0^T\mathbb{E}_{\bar{\mathbb{Q}}}\left(\sup\limits_{0\leq t\leq s}\left|(X_t-Y_t){\rm I}_{\{t\leq\tau_n\}}\right|^2\right){\rm d}s.
\eeqnn
By the Gronwall inequality, we obtain
\beqnn
\mathbb{E}_{\bar{\mathbb{Q}}}\left(\sup\limits_{0\leq t\leq T}\left|(X_t-Y_t){\rm I}_{\{t\leq\tau_n\}}\right|^2\right)
\ar\leq\ar
3(\|\sigma\|_\infty T^H)^2\exp[3K_6^2T(T^{2H}+C(2))]\mathbb{E}_{\bar{\mathbb{Q}}}\int_0^Tu_s^2{\rm d}s.\cr
\eeqnn
The Fatou lemma leads to
\beqnn
\mathbb{E}_{\bar{\mathbb{Q}}}\left(\sup\limits_{0\leq t\leq T}|X_t-Y_t|^2\right)
\leq3(\|\sigma\|_\infty T^H)^2\exp[3K_6^2T(T^{2H}+C(2))]\mathbb{E}_{\bar{\mathbb{Q}}}\int_0^Tu_s^2{\rm d}s.
\eeqnn
Hence, we deduce that
\beqnn
W_2^{d_\infty}(\mathbb{Q},\mathbb{P}_x)^2\leq2\alpha(T,H)\mathbb{H}(\mathbb{Q}|\mathbb{P}_x)
\eeqnn
with $\alpha(T,H)=3(\|\sigma\|_\infty T^H)^2\exp[3K_6^2T(T^{2H}+C(2))]$.\\
For the metric $d_2$, using the above procedure, we also can prove
\beqnn
\mathbb{E}_{\bar{\mathbb{Q}}}\left(\left|(X_t-Y_t){\rm I}_{\{t\leq\tau_n\}}\right|^2\right)
\ar\leq\ar3(\|\sigma\|_\infty T^H)^2\int_0^t\mathbb{E}_{\bar{\mathbb{Q}}}u_s^2{\rm d}s\cr
\ar+\ar3K_6^2(T^{2H}+C(2))\int_0^t\mathbb{E}_{\bar{\mathbb{Q}}}\left(\left|(X_s-Y_s){\rm I}_{\{s\leq\tau_n\}}\right|^2\right){\rm d}s.
\eeqnn
The Gronwall inequality, together with the Fatou lemma, yields
\beqnn
\mathbb{E}_{\bar{\mathbb{Q}}}\left(|X_t-Y_t|^2\right)
\leq3(\|\sigma\|_\infty T^H)^2\int_0^t\exp[3K_6^2(T^{2H}+C(2))(t-s)]\mathbb{E}_{\bar{\mathbb{Q}}}u_s^2{\rm d}s.
\eeqnn
Thus it follows that
\beqnn
W_2^{d_2}(\mathbb{Q},\mathbb{P}_x)^2\ar\leq\ar\mathbb{E}_{\bar{\mathbb{Q}}}\left(\int_0^T|X_t-Y_t|^2{\rm d}t\right)\cr
\ar\leq\ar3(\|\sigma\|_\infty T^H)^2\mathbb{E}_{\bar{\mathbb{Q}}}\int_0^Tu_s^2\left(\int_s^T\exp[3K_6^2(T^{2H}+C(2))(t-s)]{\rm d}t\right){\rm d}s\cr
\ar\leq\ar2\beta(T,H)\mathbb{H}(\mathbb{Q}|\mathbb{P}_x),
\eeqnn
where $\beta(T,H)=3(\|\sigma\|_\infty T^H)^2\frac{{\rm e}^{3K_6^2T(T^{2H}+C(2))}-1}{3K_6^2(T^{2H}+C(2))}$. The proof is complete.

\bremark\label{r5.1}
In general, $\int_0^{t\wedge\tau_n}K_H(t\wedge\tau_n,s)f(X_s,Y_s){\rm d}\bar{W}_s$ does not make sense, which forces us to consider
$\mathbb{E}_{\bar{\mathbb{Q}}}\left(\left|(X_t-Y_t){\rm I}_{\{t\leq\tau_n\}}\right|^2\right)$ rather than
$\mathbb{E}_{\bar{\mathbb{Q}}}\left(|X_{t\wedge\tau_n}-Y_{t\wedge\tau_n}|^2\right)$. Further reading on stochastic Volterra equation,
one can see \cite{Zhang10b} and references therein.
\eremark

\textbf{Acknowledgement}
The author would like to thank Professor Feng-Yu Wang for his encouragement and comments that have led to
improvements of the manuscript.
This work was supported in part by the Research project of Natural Science Foundation of
Anhui Provincial Universities (Grant No. KJ2013A134), National Natural Science Foundation of China
(Grant No. 11371029).

\end{document}